\documentclass[letterpaper,11pt]{amsart}

\usepackage{latexsym,array,delarray,amsthm,amssymb,epsfig}

\theoremstyle{plain}
\newtheorem{thm}{Theorem}[section]
\newtheorem{lemma}[thm]{Lemma}
\newtheorem{prop}[thm]{Proposition}
\newtheorem{cor}[thm]{Corollary}
\newtheorem{rem}[thm]{Remark}

\theoremstyle{definition}
\newtheorem{df}[thm]{Definition}
\newtheorem{ex}[thm]{Example}

\newtheorem{conj}[thm]{Conjecture}

\newcommand{\zz}{\mathbb{Z}}
\newcommand{\nn}{\mathbb{N}}
\newcommand{\pp}{\mathbb{P}}

\newcommand{\qq}{\mathbb{Q}}
\newcommand{\rr}{\mathbb{R}}

\newcommand{\kk}{\mathbb{K}}
\newcommand{\aaa}{\mathbf{a}}

\newcommand{\xx}{\mathbf{x}}
\newcommand{\yy}{\mathbf{y}}
\newcommand{\uu}{\mathbf{u}}
\newcommand{\vv}{\mathbf{v}}
\newcommand{\ww}{\mathbf{w}}
\newcommand{\ca}{\mathcal{A}}
\newcommand{\conv}{\mathrm{conv}}

\begin{document}

\title{Combinatorial Secant Varieties}

\author{Bernd Sturmfels}
\address{Department of Mathematics, University of California, Berkeley, CA 94720}
\email{bernd@math.berkeley.edu}
\author{Seth Sullivant}
\address{Society of Fellows and Department of Mathematics, Harvard University,
Cambridge, MA 02138}
\email{seths@math.harvard.edu}

\begin{abstract}
The construction of joins and secant varieties is studied in the combinatorial
context of monomial ideals. For ideals generated by quadratic monomials,
the generators of the secant ideals are obstructions to graph colorings,
and this leads to a commutative algebra version of the
Strong Perfect Graph Theorem. Given any projective variety and any
term order, we explore whether the 
initial ideal of the secant ideal coincides with the
secant ideal of the initial ideal. For toric varieties, this leads to the
notion of delightful triangulations of convex polytopes.
\end{abstract}

\maketitle


\section{Introduction}

Given two varieties $X$ and $ Y$
in a projective space, their {\em join} $X * Y$ is the Zariski closure of 
the union of all lines spanned by a point in $X$ and a point in $Y$.  
The join of a variety $X$ with itself is the {\em secant variety} of $X$, 
and the $r$-fold join of $X$ with itself is the {\em $r$-th secant variety} of $X$.  
It is denoted $ X^{\{r\}} = X*X* \cdots * X $. The study of
joins and secant varieties has a long tradition in algebraic geometry,
and many authors have studied the dimension and degree
of these varieties. Recent references include
\cite{Alexander1995, Catalano1996, Catalisano2002, Catalisano2005, Cox2005}.
In the emerging field of {\em algebraic statistics},
the construction of joins and secant varieties corresponds to {\em mixture models} 
\cite{Garcia2005, Pachter2005}, and it is of considerable interest to compute the defining 
prime ideals of $X*Y$ and $X^{\{r\}}$ from those of $X$ and $Y$.
For recent successes along these lines see \cite{Allman2004, Landsberg2004}.

In this paper, we present a combinatorial framework for
the study of joins and secant varieties. The basic setup
was already suggested by Simis and Ulrich \cite{Simis2000}, and our results are generalizations and extensions of theirs.
Our strategy is summarized by the following steps.  First, we take secants and joins of arbitrary projective schemes, and, hence, of arbitrary homogeneous ideals in a polynomial ring.  Second, we develop the combinatorial study of secants and joins of monomial ideals, relating secants and joins to Alexander duality, coloring properties of graphs, antichains in posets, and regular triangulations of polytopes.  Third, we use Gr\"obner degeneration as a tool to reduce questions about secants and joins of arbitrary projective schemes to secants and joins of monomial schemes.  
Among the applications of our technique is a new perspective on
classical determinantal ideals,  yielding a short unified proof for
the  Gr\"obner basis property of minors and Pfaffians.

Here is the outline for our paper.  In Section \ref{sec:joinmono} we introduce secants and joins of arbitrary ideals and we study the secants of monomial ideals.
We give an explicit formula (Theorem \ref{thm:alex1}), 
valid in characteristic zero, for
computing the join of monomial ideals 
by multiplying their Alexander duals.

 In Section \ref{sec:graph} we focus on the case of ideals
generated by quadratic monomials. If the generators are
squarefree, so the ideal is an {\em edge ideal}, then
the secants reflect coloring properties of the graph.
  As a consequence, perfect graphs make a surprise appearance,
  and we get a commutative algebra version (Theorem \ref{CASPGT})
  of the celebrated {\em strong perfect graph theorem} \cite{Chudnovsky2002}.

   In Section \ref{sec:initial} we show that the secant of an initial ideal contains 
   the initial ideal of the secant.  This allows for the derivation of numerical invariants of the secants of an ideal from the secants of carefully chosen initial ideals. We also introduce the notion of a 
   {\em delightful term order} for a variety $X$.
This is a term order where taking secants commutes with taking initial ideals. 
Diagonal term orders for determinantal and Pfaffian ideals are delightful. 

  In Section \ref{sec:triangle} we apply our techniques to the study of secant varieties of toric varieties.  We show how information about such secant varieties can be derived from regular triangulations of the corresponding polytopes. We are particularly interested
 in finding {\em delightful triangulations} which correspond to delightful term orders
 for toric varieties. The existence of delightful triangulations is explored for
 Veronese varieties, Segre varieties and scrolls.

 We close the Introduction with an example which demonstrates
 how our approach  can be used to derive equations defining secant varieties.
 Let $X \subset \pp^9$ be the cubic Veronese surface
in its standard toric embedding.
 Consider the Gr\"obner degeneration 
 of $X$ into a union of nine coordinate planes
 corresponding to the triangulation depicted in Figure 1.

\begin{figure}[h]
\begin{center}\includegraphics{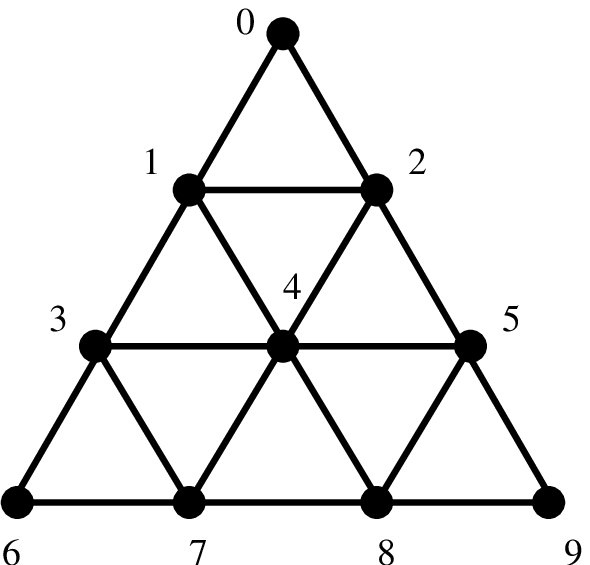}\end{center}
{\bf Figure 1}: A delightful triangulation for the cubic Veronese surface.
\end{figure}

 The initial ideal of the surface $X$ with respect to this term order
 is the edge ideal $I(G)$ whose graph $G$ consists of all
 non-edges of this triangulation:
 $$ I(G) \quad = \quad \langle
 x_0 x_3, x_0 x_4, x_0 x_5, \ldots, x_6 x_9, x_7 x_9 \rangle.
 $$
 Consider the variety  $X^{\{3\}}$
 of secant planes. This is
 a hypersurface in $\pp^9$
 and we wish to compute its defining polynomial $f$.
 To do so, we apply Theorem \ref{notcolorable} below to see that the
  ideal of the combinatorial secant variety equals
$$ I(G)^{\{3\}} \quad = \quad \langle x_0 x_4 x_6 x_9 \rangle . $$
The generator is the unique clique of size four in $G$. Equivalently, it is
the unique independent set of size four in the edge graph of the triangulation.
The desired polynomial $f$ has degree at least four, and its
leading term is a multiple of $x_0 x_4 x_6 x_9$. 
Consider the quartic invariant of ternary cubics \cite{Sturmfels1993}
\begin{eqnarray*} &
\underline{x_0 x_4 x_6 x_9} + x_1 x_2 x_7 x_8 + x_1 x_3 x_5 x_8 + x_2 x_3 x_5 x_7
- x_1^2 x_8^2 - x_2^2 x_7^2 - x_3^2 x_5^2  \\ &
-x_0 x_3 x_7 x_9 - x_0 x_4 x_7 x_8 - x_0 x_5 x_6 x_8 
- x_1 x_2 x_6 x_9 - x_1 x_3 x_4 x_9 - x_2 x_4 x_5 x_6 \\ &
+ x_0 x_3 x_8^2 + x_0 x_5 x_7^2 + x_1^2 x_7 x_9   
+ x_1 x_5^2 x_6 + x_2^2 x_6 x_8 + x_2 x_3^2 x_9 \\ &
- 3 x_1 x_4 x_5 x_7 - 3 x_2 x_3 x_4 x_8 
+ 2 x_1 x_4^2 x_8 + 2 x_2 x_4^2 x_7 + 2 x_3 x_4^2 x_5 - x_4^4.
\end{eqnarray*}
This polynomial vanishes on $X^{\{3\}}$
and it has the correct leading term.
This proves that the desired generator $f$ equals the quartic polynomial
above.

\subsection*{Acknowledgments}
The idea of using Gr\"obner degenerations as a tool to study the dimensions of secant varieties was suggested to us by Rick Miranda.   Bernd Sturmfels was supported by the NSF 
(DMS-0456960). This work was carried out while he visited
ETH Z\"urich in the summer of 2005.
Seth Sullivant was supported by an NSF graduate research fellowship.


\section{Joins of monomial ideals}\label{sec:joinmono}

Let $I_1 ,I_2, \ldots, I_r$ be ideals in a polynomial ring
$\kk[\xx] = \kk[x_1,\ldots,x_n]$ over a field $\kk$.
Their {\em join} $\,I_1 * I_2 * \cdots * I_r \,$ is a new ideal
in $\kk[\xx]$ which can be computed as follows. We introduce
$rn$ new unknowns, grouped in $r$ vectors
${\bf y}_j = (y_{j1}, \ldots,y_{jn})$, $j = 1,2,\ldots,r$, and we consider
the polynomial ring $\kk[\xx,\yy]$ in all $rn+n$ unknowns.
Let $I_j(\yy_j)$ be the image of the ideal $I_j$ in
$\kk[\xx,\yy]$ under the map $\xx \mapsto \yy_j$. Then
$I_1 * I_2 * \cdots * I_r  $ is the elimination ideal
\begin{equation}
\label{joindef}
\biggl( I_1(\yy_1) + \cdots + 
 I_r(\yy_r) + \langle \, y_{1i} + y_{2i} + \cdots + y_{ri} - x_i
 \, :  \, 1 \leq i \leq n \,  \rangle\biggr) \,\cap \, \kk[\xx]. 
\end{equation}
Of particular interest is the case when all
$r$ ideals are identical. We define the 
{\em $r$th secant} of an ideal $I \subset \kk[\xx]$ to be
the $r$-fold join of $I$ with itself:
\begin{equation}
\label{secantdef}
I^{\{r\}} \quad := \quad I * I * \cdots * I . 
\end{equation}

The join operation $I * J$ of ideals is commutative and associative.
Moreover, it satisfies the following distributive law with respect to
intersection.

\begin{lemma}
\label{distributive}
If $I,J$ and $K$ are ideals in $\kk[\xx]$ then
$$ (I  \,\cap \,J) \, *  \, K \quad = \quad
(I  * K)  \,\cap \, (J * K). $$
\end{lemma}

\begin{proof}
See Proposition 1.2 (i) in \cite{Simis2000}.
\end{proof}

If the ideals $I_1 , \ldots, I_r $ are geometrically prime  then their join $I_1 * \cdots * I_r$
is also geometrically prime. Similarly, for the properties of being 
geometrically primary and for being radical,
provided $\kk$ is a perfect field. See \cite[Proposition 1.2]{Simis2000}.
Thus  for homogeneous prime ideals, the ideal-theoretic join and secant represent the
prime ideals of the secant varieties and joins of irreducible
projective varieties, the setting discussed in the Introduction.
For arbitrary ideals, the ideal-theoretic join corresponds to the Minkowski sum of affine schemes.

This section is concerned with 
another extreme case, namely, when
the given ideals are monomial ideals. We start with
the simplest example.

\begin{ex} \label{oneVariable}
Let $n = 1$ and consider the ideals
$I = \langle x^i \rangle $ and $J = \langle x^j \rangle $.
Then $\,I * J \, = \, \langle x^k \rangle \,$
where $k $ is the smallest integer such that the
characteristic of $\kk$ divides 
${k \choose l}$ for all $l \in \{k-j+1,k-j+2,\ldots,i-1\}$.
In particular,
$$\langle x^i \rangle * \langle x^j \rangle 
\quad = \quad \langle x^{i+j-1} \rangle
\qquad \hbox{if} \ \ {\rm char}(\kk) = 0 . $$
\end{ex}

 This example generalizes to 
 {\em irreducible} monomial ideals in $n$ variables.
 Such an ideal is represented by  an integer vector
 $\uu = (u_1,\ldots,u_n) $ as follows:
 $$ {\bf m}^{\uu} \quad = \quad  \langle \,x_i^{u_i} \,: \, u_i > 0 \,\rangle.  $$

 \begin{lemma} 
 \label{joinirred}
 The join of two irreducible monomial ideals 
 ${\bf m}^\uu$ and ${\bf m}^\vv$ is
  an irreducible monomial ideal $ {\bf m}^\ww$.
  Here $w_i  = 0 $ if $u_i = 0$ or $v_i = 0$.
  Otherwise $w_i$ is the smallest integer such that the
  characteristic of   $\kk$ divides 
${w_i \choose l}$ for all $l$
with $w_i-u_i < l < v_i $, and
if ${\rm char}(\kk)=0$ then $w_i = u_i+v_i -1 $.
 \end{lemma}

 \begin{proof}
 A polynomial $f(\xx) = \sum_\aaa c_\aaa \xx^\aaa$ lies in
 the join ${\bf m}^\uu * {\bf m}^\vv$ if and only if
 $\,f(\yy_1 + \yy_2) = \sum_\aaa c_\aaa (\yy_1 + \yy_2)^\aaa \,$
lies in  the monomial ideal
 \begin{equation}
 \label{mmonoid}
 {\bf m}^\uu(\yy_1) + {\bf m}^\vv(\yy_2)
 \quad = \quad
 \langle y_{1i}^{u_i}\,:\, u_i > 0 \rangle + 
 \langle y_{2i}^{v_i}\,:\, v_i > 0 \rangle .
 \end{equation}
 This happens if and only if  every term 
$ (\yy_1 + \yy_2)^\aaa$ lies in  (\ref{mmonoid}).
Hence ${\bf m}^\uu * {\bf m}^\vv$ is  a monomial ideal.
Now, $\xx^\aaa$ lies in ${\bf m}^\uu * {\bf m}^\vv$
if and only if every term of
$$ (\yy_1 + \yy_2)^\aaa \quad = \quad
 \prod_{i=1}^n(y_{1i} + y_{2i})^{a_i}
\quad = \quad \prod_{i=1}^n \biggl(
\sum_{l=0}^{a_i} {a_i \choose l} y_{1i}^l y_{2i}^{a_i-l} \biggr) $$
lies in $\,{\bf m}^\uu(\yy_1) + {\bf m}^\vv(\yy_2) \,$
if and only if $\,w_i \leq a_i$ for some $i$ with $u_i v_i \not=0$
if and only if $\xx^\aaa $ lies in ${\bf m}^\ww$.
Therefore, $\,{\bf m}^\uu * {\bf m}^\vv \, = \,{\bf m}^\ww$.
 \end{proof}

We shall prove that the join of monomial ideals
is always a monomial ideal.  Recall that the
{\em standard monomials} of a monomial ideal $J$
are the monomials in $\kk[\xx] \backslash J$, so
$J$ is characterized by its set of standard monomials.

\begin{prop} \label{prodofstandard}
Let $I_1,  \ldots, I_r$ be monomial ideals in $\kk[\xx]$.
Then $I_1 * \cdots * I_r$ is a monomial ideal.
If ${\rm char}(\kk) = 0$ then the standard monomials 
of $I_1 * \cdots * I_r$
 are precisely the products $m_1 \cdots m_r$
where $m_j$ is standard for $I_j$.  If $I_1, \ldots, I_r$ are squarefree, the monomial generating set of $I_1  * \cdots * I_r$ is independent of ${\rm char}(\kk)$.
\end{prop}

\begin{proof}
It suffices to consider the case $r=2$; the general statement follows by induction on $r$.
If $I_1 = {\bf m}^\uu $ and $I_2 = {\bf m}^\vv$ are irreducible ideals then
both statements follow from Lemma \ref{joinirred}. Otherwise, we decompose
 $I_1  = \bigcap_\uu {\bf m}^\uu$ and 
 $I_2  = \bigcap_\vv {\bf m}^\vv$ as intersections of irreducible monomials ideals (see \cite{Miller2005}). Using
Lemma \ref{distributive}, we then write  $I_1 * I_2$ is an intersection of joins
${\bf m}^\uu * {\bf m}^\vv$. Hence $I_1 * I_2$ is a monomial ideal, and its
set of standard monomials  is the union of the sets of standard monomials
of its irreducible components ${\bf m}^\uu * {\bf m}^\vv$.  If $I_1$ and $I_2$ are irreducible and squarefree the formula for $I_1 * I_2$  of Lemma \ref{joinirred} does not depend on ${\rm char}( \kk)$.  Since every monomial ideal in the irreducible decomposition of squarefree monomial ideals is squarefree, we deduce that the monomial generators of  $I_1 * I_2$ are independent
 of the characteristic of the field $\kk$.
\end{proof}

A statement equivalent to
Proposition \ref{prodofstandard} appears in
\cite[Proposition 3.1]{Simis2000}.  

\begin{cor} \label{cor:standard}
The $r$th secant $I^{\{r\}}$
of a monomial ideal $I$ is a monomial ideal. Every
standard monomial of $I^{\{r\}}$ is a product of $r$ standard
monomials of~$I$. If ${\rm char}(\kk) = 0$ then
 every such product  is standard for $I^{\{r\}}$.  If $I$ is squarefree, the generating set of $I^{\{r\}}$ is independent of ${\rm char}(\kk)$.
\end{cor}

These results show that  the operations
of taking joins and secants are very natural from
the point of view of  Alexander duality. Namely,
forming joins is Alexander dual to taking products
of monomial ideals, and forming secants is Alexander
dual to taking powers of monomial ideals. We make
this statement precise using the $I^{[\aaa]}$ notation.
See \cite[Chapter 5]{Miller2005} for the relevant
definitions and basic facts on Alexander duality
of monomial ideals.

\begin{thm} \label{thm:alex1}
Let $I$ and $J$ be monomial ideals in
$\kk[\xx]$, ${\rm char}(\kk) = 0$,
and let $\aaa$  be a vector in $\nn^n$ 
whose coordinates are sufficiently large. Then
\begin{equation}
\label{stardual}
 I * J \quad = \quad
\bigl(I^{[\aaa]} \cdot J^{[\aaa]}\bigr)^{[2\aaa ]} 
\quad \,\, \hbox{\rm modulo} \quad
{\bf m}^{\aaa + {\bf 1}} . 
\end{equation}
Here ${\bf 1}  = (1,1,\ldots,1)$, and
the operation \ {\rm modulo ${\bf m}^{\aaa + {\bf 1}} \,$}
removes all the monomial generators
that are divisible by $x_i^{a_i+1}$ for some $i$.
\end{thm}

\begin{proof}
First assume that the given ideals are
irreducible, say, $\, I = {\bf m}^\uu$ and $J = {\bf m}^\vv$.
Then $I^{[\aaa]}$ is the principal ideal generated by
$\,\prod_{i: u_i > 0} x_i^{a_i+1-u_i}$, and $J^{[\aaa]}$ is generated by 
$\,\prod_{i: v_i > 0} x_i^{a_i+1-v_i}$. Their product is the principal ideal
$$ I^{[\aaa] }\cdot J^{[\aaa]}
\, \,  = \, \,
\biggl< \,
\prod_{i: u_i>0,v_i>0} \!\!\! x_i^{2a_i + 2-u_i-v_i}
\,\,\cdot \!\!\! \prod_{i: u_i>0,v_i=0} \!\!\! x_i^{a_i+1-u_i}
\,\,\cdot \!\!\! \prod_{i: u_i=0,v_i>0} \!\!\! x_i^{a_i+1-v_i} \,\,
\biggr>. $$ 
Taking the Alexander dual again, we see that 
$ (I^{[\aaa] }\cdot J^{[\aaa]})^{[2\aaa ]} $ is an
irreducible ideal which is generated by three
groups of monomials. The first group is
$$
x_i^{2a_i+1-(2a_i+2-u_i-v_i)}=
x_i^{u_i+v_i-1}
\qquad \hbox{for} \,\,i \,\,\hbox{such that} \,\,u_i>0 \,\,\hbox{and} \,\,v_i > 0.
$$
The second group of generators of 
$ (I^{[\aaa] }\cdot J^{[\aaa})^{[2\aaa ]} $  is
$$ x_i^{2a_i+1-(a_i+1-u_i)} \,=\, x_i^{a_i+u_i}
\qquad \hbox{for} \,\,i \,\,\hbox{such that} \,\,u_i>0 \,\,\hbox{and} \,\,v_i = 0, $$
and the third group of generators is
$$ x_i^{2a_i+1-(a_i+1-v_i)} \, = \, x_i^{a_i+v_i}
\qquad \hbox{for} \,\,i \,\,\hbox{such that} \,\,u_i=0 \,\,\hbox{and} \,\,v_i > 0. $$
Reduction modulo 
 $\,{\bf m}^{\aaa+ {\bf 1}} \,$ removes
 the second and third group of generators.
 The remaining first group generates the
 irreducible ideal $I * J$, by Lemma \ref{joinirred}.
This proves Theorem   \ref{thm:alex1} for
irreducible monomial ideals.

For the general case, we decompose the two given monomial
ideals into their irreducible components:
$\, I = \cap_\nu I_\nu \ \,$ and $\, J = \cap_\mu J_\mu $.
Alexander duality switches intersections of monomial
ideals with sum of monomial ideals, so we get
$\,I^{[\aaa]} = \sum_\nu (I_\nu)^{[\aaa]}\,$ and
$\,J^{[\aaa]} = \sum_\mu (J_\mu)^{[\aaa]}$. This implies
$$ I^{[\aaa]} \cdot J^{[\aaa]} \,=\,
\sum_{\nu,\mu}  (I_\nu)^{[\aaa]} \cdot (J_\mu)^{[\aaa]}, $$
and therefore
$$ (I^{[\aaa]} \cdot J^{[\aaa]})^{[2\aaa]} \quad = \quad
\bigcap_{\nu,\mu}  \bigl( (I_\nu)^{[\aaa]} \cdot (J_\mu)^{[\aaa]} \bigr)^{[2\aaa]}. $$
Using Lemma \ref{distributive}, and using
the result for irreducible ideals, we find
$$ I * J \,\,\, = \,\,\, \bigcap_{\nu,\mu}  (I_\nu * J_\mu )  
\,\,\, = \,\,\,  (I^{[\aaa]} \cdot J^{[\aaa]})^{[2\aaa]} \quad 
{\rm modulo}  \,\,\,\, {\bf m}^{\aaa+ {\bf 1}} . $$
This completes the proof of  Theorem \ref{thm:alex1}.
\end{proof}

\begin{cor}\label{cor:alex2}
Let $I$ be a monomial ideal in $\kk[\xx]$, suppose that ${\rm char}(\kk) = 0$,
and let ${\bf a}$ be a vector in $\nn^n$ whose coordinates are sufficiently large.
$$
I^{\{r\}} \quad = \quad  \left( ( I^{[\aaa]} )^r \right)^{[r \aaa ]} \quad \mbox{ modulo } \quad 
{\bf m}^{\aaa +  {\bf 1}}.
$$
\end{cor}

Theorem \ref{thm:alex1} and Corollary \ref{cor:alex2} can be used for the efficient computation 
of joins and secants of monomial ideals in characteristic zero.  

\begin{ex} \label{M2}
We present some code in the computer algebra program
{\tt Macaulay 2}   \cite{Grayson} for computing the
 first secant of a monomial ideal. In our example, the input is the ideal 
$\,I = \langle\, x^3, \,x^2 y^2 ,\, x z^3 ,\, y^4 ,\, y^2 z \,\rangle\,$
in $\qq[x,y,z]$. 
\begin{verbatim}
R = QQ[x,y,z]; a = 7;  
I = monomialIdeal ( x^3 , x^2*y^2 , x*z^3 , y^4 , y^2*z );
Ma  = monomialIdeal(apply(gens R, u -> u^(  a+1))); 
M2a = monomialIdeal(apply(gens R, u -> u^(2*a+1)));
Ia = monomialIdeal ((gens (Ma:I)) % Ma);   -- Alexander dualize 
Ia2 = Ia*Ia;                       -- Take square of the result
Ia22a = monomialIdeal((gens(M2a:Ia2))%M2a);-- Alexander dualize 
monomialIdeal ((gens Ia22a) % Ma)      -- reduce modulo m^{a+1}
\end{verbatim}

The output of these commands is the join of $I$ with itself:
$$\,I^{\{2\}} \,\,\,=\,\,\, I * I  \,\,\,=\,\,\,
\langle  \, x^5, \, x^4 y^3, \,x^3 y^5,\, y^7, \,y^5 z, \,x^2 y^3 z^3, \,x^3z^5
\, \rangle . $$
Note that we compute the Alexander dual in {\tt Macaulay 2} using the formula
$$\,I^{[\aaa]} \quad = \quad
({\bf m}^{\aaa+{\bf 1}} : I) 
\,\,{\rm modulo} \,\,
{\bf m}^{\aaa+{\bf 1}}. $$
\vskip -.6cm \qed
\end{ex}

The proof of Theorem \ref{thm:alex1} shows that the smallest possible choice for $\aaa$  has $a_i = \max ( 2d_i - 1, 1)$ where $d_i$ is the largest power of $x_i$ appearing in any minimal generator of $I$ or $J$.  This guarantees that none of the generators of the form $x^{u_i + v_i - 1}$ are removed when reducing modulo 
 $\,{\bf m}^{\aaa+ {\bf 1}}$.   For the secant ideal $I^{\{r\}}$, the smallest possible choice for $\aaa$ has $a_i = \max (rd_i - r + 1, 1)$ where $d_i$ is the largest power of $x_i$ appearing in any minimal generator of $I$.
 In particular, if $I$ and $J$ are squarefree monomial ideals
 we may choose $\aaa = {\bf 1}$.
Note that, for $I$ squarefree, the ideal
$I^{[{\bf 1}]} $ coincides with the squarefree Alexander dual $I^\vee$,
which is familiar from the study of Stanley-Reisner ideals $I$.
 The  code above was used with ${\tt a} = {\tt 1}$ for 
 many examples of squarefree ideals $I$
which we computed for the research presented
in the next sections.

\begin{rem} \label{complexR}
Let $\Delta$ be the simplicial complex
of $I$ and $\Delta^{\{r\}}$ the simplicial
complex of $I^{\{r\}}$.
The simplices in $\Delta^{\{r\}}$
are the unions of $r$ simplices in $\Delta$.
\end{rem}

\begin{proof}
This follows from Corollary \ref{cor:standard}. See also
\cite[Corollary 3.3]{Simis2000}.
\end{proof}


\section{Secants of edge ideals}\label{sec:graph}

Let $G$ be an undirected graph with vertex set
$[n] = \{1,2,\ldots,n\}$. To $G$ we associate the
 \emph{edge ideal} $I(G)$ which is generated
 by the squarefree quadratic monomials $x_i x_j$ 
 corresponding to the edges $\{i,j\} $ of $G$.
 For example, if $G$ is the five-cycle
 with edges $\,\bigl\{ \{1,2\}, \{2,3\}, 
 \{3,4\}, \{4,5\}, \{5,1\} \bigr\} $   then
 $$ I(G) \,\,\,=\,\,\, \langle x_1 x_2, x_2 x_3, x_3 x_4, x_4 x_5, x_5 x_1 \rangle . $$
The results below (or Remark \ref{complexR})
show that the secants of this ideal are
 \begin{equation}
 \label{fivecycle}
  I(G)^{\{2\}} \,  = \,\langle x_1 x_2 x_3 x_4 x_5 \rangle
 \quad \hbox{and} \quad
 I(G)^{\{r\}} \, = \, \langle 0 \rangle \,\,\,\hbox{for} \,\, r \geq 3 . 
 \end{equation}
 Edge ideals have been much studied in combinatorial commutative algebra.
 The main  emphasis has been on  expressing homological invariants of the
ideal $I(G)$ in terms of the graph $G$.
 In this section we relate coloring properties of the 
 graph $G$ to  algebraic properties of the secant
ideals $I(G)^{\{r\}}$.  

Recall that the {\em chromatic number}
$\chi(G)$ of a graph $G$ is the smallest number of colors which can be used to give a coloring of the vertices of $G$ such that no two adjacent vertices have the same color.  To
 the subset $V \subseteq [n]$ 
 we associate the monomial  $m_V = \prod_{i \in V} x_i$.  A basic first result is:

\begin{prop}\label{prop:chromatic}
The chromatic number $\chi(G)$ of a graph 
$G$ is the smallest integer $r \geq 0 $ such that
the $r$th secant ideal $I(G)^{\{r\}}$
is the zero ideal $ \left< 0 \right>$.
\end{prop}

\begin{proof}
The monomial $m_V = \prod_{i \in V} x_i$ is a standard monomial of $I(G)$ if and only if $V$ is an independent subset of the vertices of $G$.  An $r$-coloring is a partition $V_1, \ldots, V_r$ of the vertices of $G$ such that each $V_i$ is an independent subset of vertices of $G$. 
An $r$-coloring exists if and only if $\,x_1 x_2 \cdots x_n =  \prod_{i = 1}^r m_{V_i} $ 
is a standard monomial of $I(G)^{\{r\}}$ if and only if $I(G)^{\{r\}} = \left< 0 \right>$, since $I(G)^{\{r\}}$ is radical. 
\end{proof}

The proof of Proposition \ref{prop:chromatic} leads to a combinatorial 
description of the minimal
 generators of the secant ideals $I(G)^{\{r\}}$.  Given a subset $V \subseteq [n]$,
we write $G_V$ for subgraph of $G$ which is induced on the 
set of vertices $V$.

\begin{thm} \label{notcolorable}
The $r$th  secant $I(G)^{\{r\}}$
of an edge ideal $I(G)$ is generated by the 
squarefree monomials $m_V$ whose subgraph $G_V$ is not $r$-colorable:
$$I(G)^{\{r\}} \,\,= \,\, \langle \,m_V \, \,\, |\, \, \, \chi(G_V) > r \, \rangle.$$ 
The minimal generators of $I(G)^{\{r\}}$ are those monomials $m_V$ such that $G_V$ is not $r$ colorable but $G_U$ is $r$-colorable for every proper subset $U \subset V$.
\end{thm}

The minimal graphs that are not $2$-colorable are
the  cycles of odd length. This explains
the computation for the five-cycle in (\ref{fivecycle}).
The special case $r=2$ of Theorem \ref{notcolorable}
was already obtained in \cite[Proposition 5.1]{Simis2000}:

\begin{cor} \label{oddCycle}
The secant $I(G)^{\{2\}}$ is minimally generated 
by the monomials $m_V$ 
whose corresponding induced subgraph $G_V$ is a cycle of odd length.
\end{cor}

This implies that even for a monomial ideal $I$, there is no bound on the degrees of minimal generators of $I^{\{2\}}$ in terms of the degrees of the generators of $I$ alone. 
Furthermore, if $I$ is generated by squarefree
quadratic monomials then $I^{\{2\}}$ 
cannot have any minimal generators of even degree.

Since every graph on $\leq r$ vertices can be colored by $r$ colors,
the minimal generators of $I(G)^{\{r\}}$ have degree at least $r+1$.
This suggests the problem of characterizing 
the graphs $G$ that have the property that the secants $I(G)^{\{r\}}$ are generated in
degree $r+1$. Recall that a graph $G$ is {\em perfect} if
 the chromatic number $\chi(G_V)$ equals the clique number $\omega(G_V)$
for every subset $V \subseteq [n]$.
The {\em clique number} is the size of the largest complete subgraph.

\begin{prop}\label{prop:perfect}
A graph $G$ is perfect if and only if 
every non-zero secant ideal $I(G)^{\{r\}}$  is generated in degree $r+1$.
\end{prop}

\begin{proof}
Suppose $G$ is perfect and let  $m_V $ be  a minimal
generator of $I(G)^{\{r\}}$.  Then $G_V$ is not $r$-colorable,
i.e.~$r < \chi(G_V)$. Since $G$ is perfect, we have
$ \chi(G_V) = \omega(G_V)$ and hence $r < \omega(G_V)$.
This means there exists a subset $U \subseteq V$ such that $G_U$ is a complete subgraph $K_{r+1}$. Since $K_{r+1}$ is not $r$-colorable, the monomial $m_U $
is in $ I(G)^{\{r\}}$. Since $m_V$ is a minimal generator, we
conclude that $U = V$.
Hence $G_V = K_{r+1}$ and $m_V$ has degree $r+1$.

Conversely, if $G$ is not perfect then we pick a  subset $V \subseteq [n]$ such that 
$\,\chi(G_V) > \omega(G_V)$. We may assume that $V$ is minimal
with this property. Setting $r = \omega(G_V)$ we 
have $|V| > r + 1$. The monomial
 $m_V $ is in $I(G)^{\{r\}}$,  
whilst $m_U \notin I(G)^{\{r\}}$ for any proper subset $U \subset V$.  
Hence $m_V$ is a minimal generator of $I(G)^{\{r\}}$ which has
degree larger than $r+1$.
\end{proof}

\begin{ex}[Cyclic Polytopes]
Let $G = \overline{C}$ be the complement of a cycle $C$ of length $n > 3$; that is, the edges of $G$ are the edges not appearing in $C$.
The edge ideal $I(G)$ is a combinatorial model
for the elliptic normal curve in $\pp^{n-1}$, since the variety of $I(G)$ has degree $n$ and geometric genus 1. 

We claim that the secant ideal $I(G)^{\{r\}}$ is
the Stanley-Reisner ideal of the boundary complex
of the $2r$-dimensional cyclic polytope with $n$ vertices.
To see this, we must analyze the structure of the maximal simplices of the secant complex $C^{\{r\}}$.  If $2r \leq n$, each of these maximal simplices consists of $2r$ points.  A set $F \subset [n]$ of cardinality
$2r$ is a maximal face of $C^{\{r\}}$ if and only if 
$F$ is a union of $r$ pairwise disjoint pairs
of the form $\{\ell,\ell+1\}$ or $\{n,1\}$.
This condition on $F$ is equivalent to saying that,
for every pair $i \leq j$, with $i,j \notin F$, the cardinality of $\{i,i+1, \ldots, j-1,j \} \cap F$ is even.  Thus, by \emph{Gale's evenness condition} (e.g.  Theorem 0.7 in \cite{Ziegler1995}), the facets of $C^{\{r\}}$ are precisely the facets of the cyclic polytope.

Note that $I(G)^{\{r\}}$ is  generated in degree $r+1$,
unless $n = 2r+1$, in which case there is a single generator
in degree $2r+1$, since $G$ is perfect if $n$
is even and minimally imperfect if $n$ is odd.  This derivation of the cyclic polytope from the ``stick
elliptic normal curve'' was suggested to us by C.~Athanasiadis and F.~Santos.
\qed
\end{ex}

The most important development in graph theory in the past few years
has been the proof, announced in 2002, of Berge's Strong Perfect Graph Conjecture by
 Chudnovsky, Robertson,
Seymour and Thomas \cite{Chudnovsky2002}.
Their theorem characterizes perfect graphs in terms of
excluded induced subgraphs.

\begin{thm}[Strong Perfect Graph Theorem] 
 \label{thm:strong} The minimal imperfect graphs are precisely the
  odd holes and the complements of the odd holes. 
\end{thm}

A {\em hole} is a cycle of length greater than 3.
The secants to the edge ideal $I(G)$ detect the minimal 
imperfections in the graph $G$.
The Strong Perfect Graph Theorem implies the
following strong result on the degrees
of the minimal generators of the secant ideals $I(G)^{\{r\}}$.

\begin{cor}\label{cor:perfectgens}
Let $G$ be an imperfect graph.  Then either
\begin{enumerate}
\item $I(G)^{\{2\}}$ has a minimal generator of odd degree bigger than three, or
\item for some $r > 2$,
$I(G)^{\{r\}}$ has its minimal generators in degrees $r+1$ and  $2r + 1$ only, and 
 $I(G)^{\{s\}}$ is generated in degree $s+1$ for $s < r$.
\end{enumerate}
\end{cor}

\begin{proof}
Let $G$ be an imperfect graph. If $G$ contains an odd cycle of length $d \geq 5$
then $I(G)^{\{2\}}$ has a minimal generator of degree $d$
by Corollary \ref{oddCycle}. If $G$ contains no such odd cycle then,
by Theorem \ref{thm:strong}, the graph $G$ contains the
complement of an odd hole. Let $2r+1$ be the minimal length
of such a hole. That subgraph has chromatic number $r+1$,
so it contributes a minimal generator of degree $2r+1$ to the ideal
$I(G)^{\{r\}}$.  Theorem \ref{thm:strong} also ensures that
$I(G)^{\{r\}}$ has no generators of degree other  than $r+1$ or $2r+1$.
\end{proof}

\begin{rem} \rm  
Corollary \ref{cor:perfectgens} is, in fact,
equivalent to the Strong Perfect Graph Theorem.  
If $I(G)^{\{2 \}}$ has a minimal generator $m_V$ of odd
degree greater than $3$, the induced subgraph $G_V$ must be an odd
hole.  If $I(G)^{\{r\}}$ has a minimal generator $m_V$ of degree  $2r + 1$
and each of  $I(G)^{\{s\}}$ is generated in degree $s+1$ for $s < r$,
then the induced subgraph $G_V$ must be minimally imperfect, with
clique number $\omega(G_V) = r$ and $2r +1$ vertices.  A theorem of
Lov\'{a}sz \cite{Lovasz1972} implies that the complementary graph $\overline{G_V}$
is also minimally imperfect with clique number
$\omega(\overline{G_V}) = 2$.  Thus, $G_V$ must be the complement of an
odd hole. \qed  
\end{rem}

One family of perfect graphs 
are the incomparability graphs of posets \cite{Bollobas1998}.
If $P$ is a poset on the set $[n]$ then the edge ideal
of its incomparability graph is the {\em Stanley-Reisner ideal} of $P$
 which is defined as follows:
 $$ J(P) \quad = \quad \left<\, x_i x_j \, \, | \, \, \,{\rm neither} \,\, 
 i \leq j \,\,{\rm nor}\,\, i \geq j \,\,{\rm in}\,\, P \,\right>. $$
In words, the ideal $J(P)$ is generated by the
$2$-element antichains of $P$. That
the incomparability graph is perfect follows from
{\em Dilworth's Theorem} which states that
the size of the largest antichain of any poset $P$
equals the minimal number of chains needed to partition $P$.
Proposition \ref{prop:perfect} implies

\begin{cor} \label{cor:antichain}
Let $P$ be a poset.  Then 
any non-zero secant ideal $J(P)^{\{r\}}$
of the Stanley-Reisner ideal $J(P)$ is generated in degree $r+1$. More precisely,
$$J(P)^{\{r\}} \quad = \quad \left< \,m_A \, \, | \, \, \, A \mbox{ is an antichain  of cardinality } r+1 
\mbox{ in } P \, \right>.$$ 
\end{cor}

The secant ideals of graph ideals have other important connections to geometric constructions in the theory of graph coloring.

\begin{rem}[The Combinatorial Space of Explanations] \rm 
Given any projective scheme $X$, there is a natural rational map $\phi_r$ from the $r$-fold
free join of $X$ to the secant variety $X^{\{r\}}$.  In the statistics literature,  the (nonnegative real) preimage of a point $x$ on (the nonnegative real part of) $X^{\{r\}}$ is known
as the {\em space of explanations} for the point $x$.   See \cite{Mond2003}.

In the situation where $X = V(I(G))$ is the simplicial complex of independent sets in a graph $G$, the space of explanations has a very nice combinatorial interpretation.  Namely, if  $X^{\{r\}} = \pp^{n-1}$ and $x$ is generic then the space of explanations is a geometric realization of  ${\tt Hom}(G,K_r)$, the polyhedral cell complex of graph homomorphisms 
from $G$ to $K_r$.  See \cite{Babson2003} and \cite[\S 4.1]{Kozlov2005}. \qed
\end{rem}

The graph-theoretic interpretation of secant ideals extends to arbitrary squarefree
monomial ideals, by thinking of these as facet ideals as in \cite{Faridi2002}.
 Let $H \subset 2^{[n]}$ be a collection of subsets
of $[n]$ with the property that for every $U,V \in H$, neither $U
\subset V$ nor $V \subset U$.  The collection $H$ is the set of hyperedges
of a hypergraph, or the maximal faces of a simplicial complex. 
The {\em facet ideal} of the hypergraph $H$ is
  the squarefree monomial ideal
$$I(H)  \, \, = \, \, \left< \,m_V \,\, |\, \, V \in H \,\right>.$$
 A coloring of the
vertices of the hypergraph is an assignment of colors with the property
that no hyperedge has all its vertices the same color.  The chromatic
number of $H$ is the smallest number $\chi(H)$ such that $H$ has a
coloring using $\chi(H)$ colors. 
Proposition \ref{prop:chromatic} easily generalizes to this setting:

\begin{prop}  
The chromatic number of a hypergraph $H$ is the smallest
positive integer $r $ such
that $I(H)^{\{r\}} = \left< 0 \right>$.
\end{prop}

In the next section we shall apply these results
to the study of secant varieties of certain irreducible
projective varieties. It might make sense to use
graph-theoretic language even at that level of generality.
We could say that a projective scheme $X$ is
{\em perfect} if all its secant schemes
are cut out by equations of minimal degree,
and the smallest secant scheme of $X$
which fills the ambient projective space
would determine the {\em chromatic number} of $X$.
We are inclined to speculate that some version of
 the Strong Perfect Graph Theorem
 generalizes to arbitrary projective schemes defined by quadrics.
One piece of evidence is that the result about the generating degrees of secants in Corollary \ref{cor:perfectgens} generalizes to arbitrary quadratic monomial ideals.

 \begin{thm} \label{CASPGT}
 Let $I$ be an ideal generated by quadratic (not necessarily squarefree) monomials 
 whose projective scheme  is not perfect and let ${\rm char}(\kk) = 0$.  Then either
\begin{enumerate}
\item $I^{\{2\}}$ has a minimal generator of odd degree bigger than three, or
\item for some $r > 2$, the ideal
$I^{\{r\}}$ has its minimal generators in degrees $r+1$ and $2r + 1$ only, and
 $I^{\{s\}}$ is generated in degree $s+1$ for $s < r$.
\end{enumerate}
\end{thm}

\begin{proof}
Fix a large integer $m \gg 0$ and introduce a graph $G_{m}$
as follows. There is one vertex $X_{i,0}$ for each index $i$
such that $x_i^2 \not\in I$, and there are $m$ vertices
$X_{i,1}, \ldots,X_{i,m}$ for each index $i$ such that $x_i^2 \in I$.
Two distinct vertices $X_{i,j}$ and $X_{i',j'}$ are connected
by an edge in $G_m$ if and only if $x_i x_{i'} \in I$.
We consider the edge ideal $I(G_m)$ in the polynomial
ring with variables $X_{i,j}$.

We claim that for every integer $r \leq m$, the ideal $I^{\{r\}}$ is obtained
from the squarefree ideal $I(G_{m})^{\{r\}}$ by replacing $X_{i,j}$ by $x_i$.
Since Theorem \ref{CASPGT} holds for 
$I(G_m)$, we conclude that it holds for the given ideal $I$ as well.

To prove the claim, note that $x_i^{r+1}$ divides a minimal generator of $I^{\{r\}}$ if and only if $x_i^{r+1}$ is a minimal generator of $I^{\{r\}}$ if and only if $x_i^2 \in I$.  This chain of implications follows because $I$ is generated by quadrics and either $x_i^2 \notin I$, in which case $I^{\{r\}}$ is squarefree in $x_i$, or $x_i^2 \in I$ in which case $x_i^{r+1} \in I^{\{r\}}$.  Since $I(G_m)^{\{r\}}$ contains a generator that maps onto $x_i^{r+1}$ if and only if $x_i^2 \in I$, it suffices to show that the replacement procedure sends $I(G_m)^{\{r\}}$ to $I^{\{r\}}$ when restricted to those monomials not divisible  by $x_i^{r+1}$ for any $I$.  Let $M $ be such a monomial and suppose
that $M$ is not in $I^{\{r\}}$ and
 that $M$ is squarefree in each variable $x_i$ such that $x_i^2 \notin I$. 
This condition holds if and only if $M$ admits a factorization
  $M = M_1 \cdots M_r$ into a product of $r$ monomials such that each $M_j$ is standard for $I$.  By our assumption on $M$, each $M_j$ is squarefree.  Since the squarefree standard monomials of $I(G_m)$ map onto the set of squarefree standard monomials of $I$, we deduce that $M \notin I^{\{r\}}$ if and only if every squarefree preimage of $M$ is standard for $I(G_m)^{\{r\}}$.  The existence of such squarefree preimages is guaranteed because $r \leq m$.  
We conclude that every minimal generator of $I^{\{r\}}$ arises from a
squarefree minimal generator of $I(G_m)^{\{r\}}$ and hence,
by Corollary \ref{cor:perfectgens},  satisfies
the specified requirements on its degree.
\end{proof} 

Needless to say, it would be fantastic to find a  commutative
algebra proof of Theorem \ref{CASPGT} and hence of the
Strong Perfect Graph Theorem. Note, however, that the statement
of Theorem \ref{CASPGT} does not hold for
all ideals generated  by quadrics.  In particular, for a 
non-monomial ideal $I$ generated by quadrics, the secant 
ideal $I^{\{2\}}$ need not have
a generator of odd degree.

\begin{ex}
Let $I$ be generated by two generic homogeneous quadrics in $ \kk[x,y,z]$.
The variety of $I$ consists of four points in general position in $\pp^2$.
The secant variety is the reduced union of six lines. Hence
  $I^{\{2\}} $ is a principal ideal generated by
  a homogeneous polynomial of degree six. \qed
  \end{ex}


\section{Equations of secant varieties via initial ideals}\label{sec:initial}

In this section we consider the degeneration of secant ideals to their initial ideals.  
Theorem \ref{thm:inclusion} and its corollaries on initial degrees 
appear already in \cite{Simis2000}, but we include some proofs because 
ours are, perhaps, more elementary. We then
examine determinantal and Pfaffian ideals. In contrast to the
exposition in \cite[\S 5]{Simis2000}, we offer 
direct new proofs for the Gr\"obner basis properties
of determinants and Pfaffians, using results from Section 3.
To be precise, we replace the use of
the Knuth-Robinson-Schensted correspondence, first proposed in \cite{Sturmfels1989},
by Dilworth's Theorem (Corollary \ref{cor:antichain}).
Besides Dilworth's Theorem, which is a relatively easy
combinatorial result,
our derivation depends only on elementary linear algebra.
The full proof for generic matrices is presented in Theorem 
\ref{NoMoreKRS} and Corollary  \ref{cor:diagonalGB} below.

Let $I_1,\ldots,I_r$ be arbitrary ideals in $\kk[\xx]$  and $\prec$ any
term order.  Then the initial ideal of a join is contained in the join of
the initial ideals. 

\begin{thm}
\label{thm:inclusion}
 We have the following inclusions of monomial ideals:
$$  {\rm in}_\prec \bigl( I_1 * I_2 * \cdots * I_r \bigr)
\quad \subseteq \quad
{\rm in}_\prec(I_1) *  {\rm in}_\prec(I_2) * \cdots * 
{\rm in}_\prec(I_r)  .$$
\end{thm}

\begin{proof}
It suffices to consider the case of the join of two ideals $I * J$; the general result following by induction on $r$. Consider any polynomial $f \in I*J$.
 Let $w \in \rr^n$ be a weight vector which represents the
term order $\prec$ in the sense that
${\rm in}_w(I) = {\rm in}_\prec(I)$, 
${\rm in}_w(J) = {\rm in}_\prec(J)$ and
${\rm in}_w(f) = {\rm in}_\prec(f)$. We 
denote the latter monomial by $\, m = {\rm in}_w(f)$.
We consider the ideal $I(\xx) + J(\yy)$ in the
 polynomial ring $\kk[\xx,\yy]$. The $(w,w)$-initial ideal of
 this ideal equals
 \begin{equation}
\label{threeideals}
{\rm in}_{(w,w)}(I(\xx) + J(\yy)) \, = \,
{\rm in}_w(I(\xx)) + {\rm in}_w(J(\yy)) \, = \,
 {\rm in}_\prec(I(\xx)) + {\rm in}_\prec(J(\yy)). 
 \end{equation}
This is seen by refining $(w,w)$ to a 
term order and using Buchberger's First Criterion
(the S-pairs of polynomials with relatively prime
leading terms reduce to zero).
Now, since $f \in I*J$, the
polynomial $f(\xx + \yy)$ lies in
$I(\xx) + J(\yy)$. Its $(w,w)$-leading form equals
$\,{\rm in}_{(w,w)}\bigl( f(\xx + \yy) \bigr) \, = \,  m(\xx+\yy)$.
This polynomial lies in (\ref{threeideals}) and hence
$m $ lies in $\,{\rm in}_\prec(I) * {\rm in}_\prec(J)$, as desired.
\end{proof}

In the special case when all $r$ ideals are equal, this proposition implies

\begin{cor}\label{cor:secin} A secant of an initial ideal contains the
initial ideal of the corresponding secant ideal. For any ideal $I$,
term order $\prec$ and integer  $r \geq 2$,
\begin{equation}
\label{secantinclusion}
{\rm in}_\prec( I^{\{r\}} ) \quad \subseteq \quad 
 ({\rm in}_\prec(I))^{\{r\}} 
\end{equation}
\end{cor}

 Theorem \ref{thm:inclusion} implies lower bounds on the degrees of generators
for the ideals of joins and secants of arbitrary projective schemes. 
 For a homogeneous ideal $I \subset \kk[\xx]$ let
${\rm indeg}(I)$ denote the smallest degree of any minimal generator of $I$.  We omit
  the proofs which appear in \cite[Theorem 4.4]{Simis2000}. 

\begin{cor} \label{cor:joinBound}
Let ${\rm char}(\kk) = 0$ and let $I$ and $J$ be homogeneous ideals
in $\kk[\xx]$.
Then either $I*J = \left< 0 \right>$ 
or  $\,{\rm indeg}(I*J) \geq  {\rm indeg}(I) + {\rm indeg}(J) - 1$.
\end{cor}

\begin{cor}\label{cor:secdegrees}
Let ${\rm char}(\kk) = 0$ and $I$ be a homogeneous ideal such that
 ${\rm indeg}(I) = d$.  Then either
$I^{\{r\}} = \left< 0 \right>$ or ${\rm indeg}(I^{\{r\}}) \geq rd - r +1 $.
\end{cor}

\begin{rem} \rm
The lower bound of Corollary \ref{cor:secdegrees} is best
possible. This is illustrated by the family of determinantal ideals
to be featured below. Namely, if $I$ is the ideal
of $d \times d$-minors of an $m \times m$-matrix of unknowns
(for $m \gg 0 $) then $I^{\{r\}}$ is the ideal of 
$(rd-r+1) \times (rd-r+1)$-minors. 
\end{rem}

\begin{rem} \rm
Corollaries \ref{cor:joinBound} and  \ref{cor:secdegrees}
do not hold if the field $\kk$ has positive characteristic.
For instance, take $n=1$ and ${\rm char}(\kk) = 2$, 
and consider the ideal $I = \langle x^3 \rangle$. We have
 ${\rm indeg}(I) = d = 3$.  By Example \ref{oneVariable}, the first 
secant ideal is $\,I^{\{2\}} \,=\, \langle x^4 \rangle \,$ while
the bound in Corollary \ref{cor:secdegrees}
says ${\rm indeg}(I^{\{2\}}) \geq 5 $. 
\end{rem}

Corollary \ref{cor:secin} shows that the secant of the initial ideal
$({\rm in}_\prec(I))^{\{r\}}$ can provide useful bounds on numerical
invariants of the ideal $I^{\{r\}}$.  An inclusion of monomial
ideals leads to a coefficientwise inequality among the Hilbert
series and hence  among values of the
Hilbert polynomials. This implies:

\begin{cor} \label{cor:dims}
We have the following inequality for the Krull dimension:
$$\dim \kk[\xx]/ ({\rm in}_\prec(I))^{\{r\}} \,\,\,\leq \,\,\,\dim \kk[\xx]/I^{\{r\}}. $$
If these two algebras have the same Krull dimension then their degrees satisfy
$$ \deg  \kk[\xx]/ ({\rm in}_\prec(I))^{\{r\}} \,\,\,\leq \,\,\,\deg \kk[\xx]/I^{\{r\}}.$$
\end{cor}

\begin{df}
If equality holds in (\ref{secantinclusion}) then we say that
the term order $\prec$ is {\em $r$-delightful} for the ideal $I$.
We call $\prec$ {\em delightful} for $I$ if this holds
for all integers $ r \geq 2$.
Being delightful implies that equalities hold in Corollary \ref{cor:dims}.
\end{df}

A classical result in combinatorial commutative algebra states 
  that  the $k \times k$ minors of a generic
matrix, the $k \times k$ minors of a generic symmetric
matrix, and the $2k \times 2k$ sub-Pfaffians of a generic
skew-symmetric matrix are all Gr\"obner bases for the ideals they
generate.  As a corollary, one deduces that these ideals are all prime
ideals and one gets formulas for their Hilbert series.
Our approach through secants of initial ideals provides a 
unified framework for proving these results, using the following strategy:
\begin{enumerate}
\item Solve the ``easy'' $k = 2$ case by specifying 
a quadratic Gr\"obner basis for $I$ whose leading terms 
correspond to the incomparable pairs in a poset $P$.
(Usually, one here has an {\em algebra with straightening law}.)
\item Determine a combinatorial description of the antichains
of size $r+1$ in $P$.
(By Corollary \ref{cor:antichain}, these  antichains generate
  $({\rm in}_\prec(I) )^{\{r\}}$).
 \item Find a set $\mathcal{G} \subset I^{\{r\}}$ whose initial terms 
  are the above antichains.
  \item Conclude  that $({\rm in}_\prec(I))^{\{r\}} = {\rm in}_\prec(
  I^{\{r\}} )$ and $\mathcal{G}$ is a Gr\"obner basis.
\end{enumerate}

The following theorem was first proved in \cite{Sturmfels1989}
using the Knuth-Robinson-Schensted correspondence.
In our new proof, Knuth-Robinson-Schensted is replaced 
by Dilworth's Theorem  (imcomparability graphs are perfect).

\begin{thm} \label{NoMoreKRS}
Let $I$ be the ideal generated by the $2 \times 2$ minors of a generic $m \times n$ matrix,
and let $\prec$ be any term order on $\kk[x_{11}, \ldots, x_{mn}]$ which selects the diagonal leading term of each $2 \times 2$ minor.  Then $\prec$ is delightful for $I$.
\end{thm}

\begin{proof}
The poset $P$ is the product of an $m$-chain with an $n$-chain,
indexed so that the incomparable pairs are $x_{i j} x_{kl}$ with
$i < k$ and $j < l$. One easily checks that ${\rm in}_\prec(I) = J(P)$.
By Corollary \ref{cor:antichain}, $J(P)^{\{r\}}$ is generated
by the monomials $x_{i_0 j_0} x_{i_1 j_1} \cdots x_{i_r j_r}$
with $i_0 < i_1 < \cdots < i_r$ and
        $j_0 < j_1 < \cdots < j_r$.
Each such monomial is the $\prec$-leading term of
an $(r+1) \times (r+1)$-minor of the $m \times n$-matrix. 

The affine variety $V(I)$ consists of all matrices of rank $\leq 1$.
Since a matrix has rank $\leq r$ if and only if it is a sum
of $r$ matrices of rank $\leq 1$,  the affine variety
$V(I^{\{r\}})$ consists of all matrices of rank $\leq r$. 
Hence the $(r +1) \times (r+1)$ minors  vanish on $V(I^{\{r\}})$.
Now, the ideal $I$ is easily seen to be prime over any field,
and hence $I^{\{r\}}$ is geometrically prime. Hence
the $(r +1) \times (r+1)$ minors lie in the ideal $\,I^{\{r\}}$.
This proves that the monomial ideal
$\,J(P)^{\{r\}} \, = \, ({\rm in}_\prec(I))^{\{r\}}\,$
is equal to the monomial ideal $\, {\rm in}_\prec(I^{\{r\}})\,$ for
all $r \geq 2$. We conclude that the term order
 $\prec$ is delightful for the ideal $I$
of $2 \times 2$-minors.
\end{proof}

\begin{cor} \label{cor:diagonalGB}
The secant ideal $I^{\{k-1\}}$ is generated by the $k \times k$ 
minors of a generic matrix, and these minors are a Gr\"obner basis
under any diagonal term order $\prec$ as above.
\end{cor}

\begin{proof}
In the proof of Theorem \ref{NoMoreKRS} we have argued
that the $k \times k$-minors lie in $I^{\{k-1\}}$,
and their leading terms generate the initial ideal
$ \, ({\rm in}_\prec(I))^{\{k-1\}}\,= \, {\rm in}_\prec(I^{\{k-1\}})$.
This implies that the $k \times k$-minors form a 
Gr\"obner basis for the ideal $I^{\{k-1\}}$, and,
in particular, they generate that ideal.
\end{proof}

\begin{cor} The ideal of $k \times  k$ minors of a generic 
matrix is prime.
\end{cor}

\begin{proof} The ideal $I$ of $2 \times 2$-minors is geometrically prime.
The secant ideal $I^{\{k-1\}}$ of a geometrically prime ideal $I$ is prime. 
Now use Corollary \ref{cor:diagonalGB}.
\end{proof}

The same argument works also for symmetric minors and Pfaffians.

\begin{ex}[Minors of a symmetric matrix]
\label{ex:symmetric}
Consider a generic $m \times m$ symmetric matrix $(x_{ij})$
and let $I$ be its ideal of $2 \times 2$-minors.
Let $P$ be the poset on  the set of variables 
$\,\{x_{ij} \,\,|\,\,\, 1 \leq i \leq j \leq m \}\,$ defined by
$\,x_{ij} \leq x_{kl}\,$ whenever $i \leq k$ and $j \geq l$.
Let $\prec$ be the reverse lexicographic term order on
any linear extension of $P$. It is easy to check that 
the $2 \times 2$-minors are a Gr\"obner basis 
for $I$ with respect to $\prec$, and the
generators of ${\rm in}_\prec(I)$ are the incomparable
pairs in $P$. Every antichain of size $k$ in $P$
is the leading term of a $k \times k$-subdeterminant
of $(x_{ij})$. Hence the term order $\prec$ is
delightful for $I$, and we conclude that
the  $k \times k$-minors of $(x_{ij})$
form a Gr\"obner basis of $I^{\{k-1\}}$ with respect to $\prec$,
and their ideal is prime.
\qed \end{ex}

\begin{ex}[Pfaffians]\label{ex:pfaff}
Consider a generic $m \times m$ skew-symmetric matrix $(x_{ij})$,
and let $I$ be the ideal generated by its $4 \times 4$-Pfaffians
$\, x_{il} x_{jk} - x_{ik} x_{jl} + x_{ij} x_{kl}\,$
for $ 1 \leq i < j < k < l \leq m$. These are the
{\em three-term Pl\"ucker relations}, and $I$ is the defining ideal of the
Grassmannian of lines in projective $(m-1)$-space and, hence, is geometrically prime.
Let $P$ be the poset on  the variables 
$\,\{\,x_{ij} \,|\, 1 \leq i <  j \leq m \}\,$ defined by
$\,x_{ij} \leq x_{kl}\,$ whenever $i \leq k$ and $j \leq l$.
Let $\prec$ be the reverse lexicographic term order on
any linear extension of $P$.
The three-term Pl\"ucker relations are a Gr\"obner basis 
for $I$ with respect to $\prec$, and the
generators of ${\rm in}_\prec(I)$ are the incomparable
pairs in $P$ (see, for example, \cite[Chapter 14]{Miller2005}). Every antichain of size $k$ in $P$
is the leading term of a $2k \times 2k$ subpfaffian
of $(x_{ij})$ and each $2k \times 2k$ subpfaffian lies in $I^{\{k-1\}}$. Hence the term order $\prec$ is
delightful for $I$. We conclude that
the  $2k \times 2k$ subpfaffians of $(x_{ij})$
form a Gr\"obner basis of $I^{\{k-1\}}$ with respect to $\prec$,
and their ideal is prime.
\qed \end{ex}

\begin{rem} \rm
We do not know whether the Pl\"ucker ideals of the higher
Grassmannians, $G_{k,n}$ for $k \geq 3$, and their Schubert
subvarieties, admit delightful term orders.
Some computational explorations would be worthwhile.
\end{rem}

\begin{rem} \rm
The arguments we have presented also work to show the Gr\"obner basis property for ladder determinantal ideals and ladder Pfaffian ideals.  These ladder ideals consist of ideals generated by the minors and Pfaffians contained in staircase shaped regions of a generic matrix, symmetric matrix, or skew-symmetric matrix.  Each poset $P$ for these ideals  is a sub-poset of the posets described in Theorem \ref{NoMoreKRS} and Examples \ref{ex:symmetric} and \ref{ex:pfaff}. 
For studies of such ideals and their posets we refer to
\cite{Conca1995} and \cite{Gonciulea2000}.
\end{rem}


\section{Delightful triangulations of polytopes}\label{sec:triangle}

In this section we consider the case when 
$I$ is a homogeneous toric ideal,
and we examine when there exist delightful initial ideals for these toric ideals.  In the case where ${\rm in}_\prec(I)$ is generated by squarefree monomials, this corresponds to finding a special regular unimodular triangulation of the point configuration
underlying $I$. We begin by briefly reviewing the connection between toric initial ideals and regular triangulations \cite{Sturmfels1995}.

Let $\mathcal{A} = \{a_1, \ldots, a_n\} \subset \zz^d$ and suppose there is a
 vector $\omega \in \qq^d$ such that $\omega^T a_i = 1$ for all $i$.
 The toric ideal $I_\ca \subset \kk[\xx]$ is the kernel of the map
$$\kk[x_1, \ldots, x_n] \rightarrow \kk[t_1^{\pm 1}, \ldots, t_d^{\pm 1}]
\,,\,\, x_j \mapsto  \prod_{i = 1}^d t_i^{a_{ij}}.$$
Let $\prec$ be any term order on $\kk[\xx]$ 
and ${\rm in}_\prec(I_\ca)$ the initial ideal of $I_\ca$.
Then the radical of   ${\rm in}_\prec(I_\ca)$ is a squarefree monomial
ideal whose corresponding simplicial complex $\Delta_\prec(\ca)$
is a regular triangulation of $\ca$. Conversely, every regular triangulation
of $\ca$ has the form $\Delta_\prec(\ca)$ for some term order $\prec$.
A subset $\{a_{i_1},\ldots,a_{i_r}\}$ of $\ca$ is a simplex of 
the triangulation $\Delta_\prec(\ca)$ if and only if
every power of $x_{i_1} \cdots x_{i_r}$ is a standard monomial modulo 
${\rm in}_\prec(I_\ca)$.

A triangulation of the point configuration $\ca$ is 
said to be \emph{full} if every point of $\ca$ appears as 
the vertex of some simplex in the triangulation.

\begin{prop}
Suppose that $\prec$ is a delightful term order for 
the toric ideal $I_\ca$.  Then
the regular triangulation $\Delta_\prec(\ca)$ is full.
\end{prop}

\begin{proof}
If $\Delta_\prec(\ca)$ is not full then some $a_i$ 
is not a vertex of $\Delta_\prec(\ca)$.  Hence
$x_i^m \in {\rm in}_\prec(I_\ca)$ for some $m> 1$.
By Example \ref{oneVariable},
$({\rm in}_\prec(I_\ca))^{\{r\}}$ contains the monomial $x_i^{rm -r + 1}$.  Thus
 $\,({\rm in}_\prec(I_\ca))^{\{r\}} \, \neq \, {\rm in}_\prec(I_\ca^{\{r\}}) = \langle 0 \rangle \,$
 for $r \gg 0$.
\end{proof}

To illustrate the notion of delightful triangulations,
and to tie it in with determinantal ideals,
 we start out with examples in dimension two $(d=3)$.

\begin{ex} Consider the embedding of 
the toric surface $\pp^1 \times \pp^1$
in $\pp^8$ by the line bundle $\mathcal{O}(2,2)$.
 Here $n= 9$ and the defining configuration is 
$$ \ca \quad  = \quad
\begin{pmatrix}
1 & 1 & 1 & 1 & 1 & 1 & 1  & 1 & 1 \\
0 & 0 & 0 & 1 & 1 & 1 & 2 & 2 & 2 \\
0 & 1 & 2 & 0 & 1 & 2 & 0 & 1 & 2 
\end{pmatrix} .$$
The toric ideal $I_\ca$ is generated by the $2 \times 2$-minors
of the symmetric matrix
$$M \quad = \quad \begin{pmatrix}
x_1 & x_2 & x_4 & x_5 \\
x_2 & x_3 & x_5 & x_6 \\
x_4 & x_5 & x_7 & x_8 \\
x_5 & x_6 & x_8 & x_9 
\end{pmatrix}. $$
Fix a term order $\prec$ which selects the main diagonal
product as the leading term for each $2 \times 2$-minor.
The regular triangulation $\Delta_\prec(\ca)$ is 
displayed on the left in Figure 2.
The diagram on the right  of Figure 2 shows a poset $P$ to which
 we can apply  steps (1)--(4) of Section 4. 
Note that the maximal chains in $P$ are the triangles in 
$\Delta_\prec(\ca)$.
Using {\tt Macaulay 2} we can verify that
the $r \times r$-minors of $M$
 form a Gr\"obner basis of $I_\ca^{\{r\}}$.
In particular, the variety of secant planes 
to our surface
is the hypersurface $\,{\rm det}(M) = 0$. 
This proves that this triangulation of $\ca$ is
delightful. Note that this example is a specialization of
the $m=4$ case in Example \ref{ex:symmetric}. 
\qed \end{ex}

\begin{figure}[ht] \label{fig:square}
\begin{center}
\includegraphics{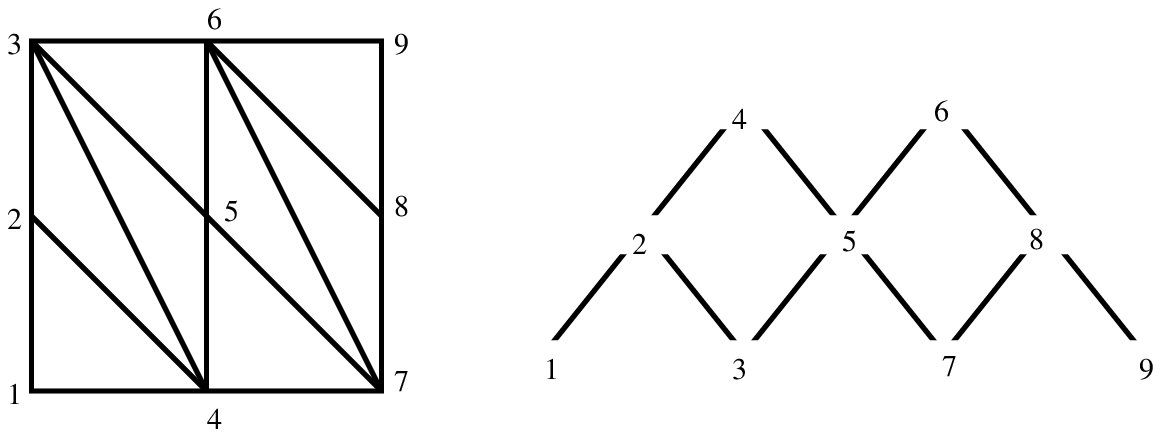}
\end{center}
{\bf Figure 2.} A delightful triangulation given by the
chains in a poset.
\end{figure}

\begin{ex}
The Veronese example discussed in the Introduction has
$$ \ca \quad = \quad
\begin{pmatrix}
3 & 2 & 2 & 1 & 1 & 1 & 0 & 0 & 0 & 0 \\ 
0 & 1 & 0 & 2 & 1 & 0 & 3 & 2 & 1 & 0 \\ 
0 & 0 & 1 & 0 & 1 & 2 & 0 & 1 & 2 & 3 
\end{pmatrix}.
$$
We have seen that the standard triangulation of $\ca$ is $3$-delightful.
However, it turns out that no full triangulation of $\ca$ is $2$-delightful.
This can be proved by a brute-force enumeration of all full triangulations
of $\ca$, using CaTS \cite{Jensen2003} or TOPCOM \cite{Rambau2002},
and by using the following counting argument.

Results in \cite{Cox2005} imply that $\,I_\ca^{\{2\}}\,$ has
the expected dimension (Krull dimension $6$) and
its degree equals $15$.
For each of the triangulations $\Delta_\prec(\ca)$ we count
the number of six-tuples of vertices which form the vertices
of two disjoint triangles. If there is no such six-tuple then
${\rm in}_\prec(I_\ca)^{\{2\}}$ has dimension less than six,
so $\prec$ cannot be delightful. Otherwise, the number of
such six-tuples equals the degree of 
${\rm in}_\prec(I_\ca)^{\{2\}}$. Now, the maximum 
number arising from any triangulation of $\ca$
is $14$ which is less than $15$.
\qed \end{ex}
A familiar example of a delightful triangulation is the
staircase triangulation of the product of two simplices.
This is the content of Theorem \ref{NoMoreKRS}.

The  quest for delightful triangulations is a worthwhile
undertaking even if
no such triangulation exists. Namely, the same approach
can be used for showing 
that certain secant varieties of toric varieties are nondefective and to
compute a non-trivial lower bound on their degree.
Recall that a $(d-1)$ dimensional subvariety $X$ of
$\pp^{n-1}$ is called $r$-{\em defective} if the
secant variety $X^{\{r\}}$ has
dimension less than $\min(rd - 1, n-1)$, which is
the {\em expected dimension}.
If all secant varieties $X^{\{r\}}$ have the expected
dimension then $X$ is called {\em nondefective}.    Regular triangulations
of $\ca$ can be used to prove that the toric variety $X_\ca$ is
nondefective.  This was the original problem suggested to us by 
Rick Miranda.  To make the idea precise, we introduce
the following terminology.

Let $\Delta$ be a full triangulation of a 
configuration $\ca$ of maximal rank $d$.
A subset $C$ of $\ca$ is
called {\em $r$-partitionable} if $C$ is the disjoint union
of $r$ maximal simplices in $\Delta$.
Naturally, if $C$ is  $r$-partitionable then $|C| = rd $.
We write $X_\ca$ for the projective toric variety in $\pp^{n-1}$
defined by $I_\ca$.

\begin{thm} \label{partitionable}
Let $\Delta$ be a regular triangulation of $\ca $
which has at least one  $r$-partitionable set.
Then $X_\ca^{\{r\}}$ has the expected dimension, and
the degree of $X_\ca^{\{r\}}$ is bounded below by the number of
$r$-partitionable sets in $\Delta$.
\end{thm}

\begin{proof}
The $r$-partitionable sets are the
$(rd-1)$-dimensional simplices of
$\Delta^{\{r\}}$ by Remark \ref{complexR}.
The number of $r$-partitionable sets is positive if and only if
$\Delta^{\{r\}}$ has the expected dimension $rd-1$.
In this case, that number is the number of maximal-dimensional
simplices in $\Delta^{\{r\}}$, which is the degree of 
$\Delta^{\{r\}}$ when regarded as a reduced union of
coordinate subspaces in $\pp^{n-1}$.

Pick a term order $\prec$ such that
$\Delta = \Delta_\prec(\ca)$. Then we have
\begin{equation}
\label{ineqchain}
 {\rm deg}(X_\ca^{\{r\}})  
\,\,=\,\,  {\rm deg}({\rm in}_\prec(I_\ca^{\{r\}}))
\,\, \geq\,\,   {\rm deg}({\rm in}_\prec(I_\ca)^{\{r\}})
\,\, \geq \,\, {\rm deg}(\Delta^{\{r\}}).
\end{equation}
The first equation holds because
the degree is preserved under Gr\"obner
degenerations, the middle inequality holds
by Corollary \ref{cor:dims}, and the
last inequality holds because 
$\Delta^{\{r\}}$ is the reduced scheme
defined by the (possibly non-radical)
ideal ${\rm in}_\prec(I_A)^{\{r\}}$.
This proves the asserted lower bound.
\end{proof}

\begin{conj}
If the lower bound for the degree in Theorem \ref{partitionable} holds with equality then  
$\Delta$ is an $r$-delightful triangulation of $\ca$. 
\end{conj}

\begin{ex}[Segre varieties, lex triangulations, and rook placements]
Let ${\bf d} =  (d_1, \ldots, d_n) $ be a vector of positive integers
and fix the configuration
$$\mathcal{A}_{\bf d} = \{v_{i_1\cdots i_n} =  {\bf e}_{i_1} \oplus   \cdots \oplus {\bf e}_{i_n} \, \, 
| \, \, 0 \leq i_j \leq d_j \mbox{ for all } j \}.$$
Thus $\mathcal{A}_{\bf d}$ represents a product of simplices, and
the corresponding toric variety is the product $\pp^{d_1} \times \pp^{d_2}  \times \cdots \times \pp^{d_n}$ in the standard Segre embedding.

Consider a lexicographic term order $\prec$ such that  $v_{i_1 \cdots i_n}$ is
higher than all other elements of $\mathcal{A}_{\bf d}$.  Since the polytope 
${\rm conv}( \mathcal{A}_{\bf d})$ is smooth, the resulting 
{\em lexicographic triangulation}
$\Delta_\prec( \mathcal{A}_{\bf d})$
  has exactly one maximal simplex which contains the vertex $v_{i_1 \cdots i_n}$.  This simplex
  is denoted $\sigma_{i_1 \cdots i_n}$, and it is
   formed by the vertices that are neighbors of $v_{i_1 \cdots i_n}$.
    In other words, the
  simplex $\sigma_{i_1 \cdots i_n}$
   contains all $v_{j_1 \cdots j_n}$ such that the Hamming distance between the vectors $(i_1, \ldots, i_n)$ and $(j_1, \ldots, j_n)$ is at most one.

Now consider a set of indices $I = \{{\bf i}_1, \ldots, {\bf i}_s \}$
with the property that the Hamming distance between ${\bf i}_j$
and ${\bf i}_k$ is greater than two for all $j \not= k$.  Let $\Delta$ be
any lexicographic triangulation of $\mathcal{A}_{\bf d}$ which puts
the elements in $V_I = \{v_{\bf i}\, |\, {\bf i} \in I \}$
lexicographically larger than all elements of
$\,\mathcal{A}_{\bf d} \backslash V_I$.   
By our assumption on the Hamming distance between elements of $I$,
each simplex $\sigma_{\bf i}, \, {\bf i} \in I \,$ appears in the
triangulation $\Delta$, and these
simplices are disjoint.  Thus, if such an index set of cardinality
$s$ exists, the secant varieties $X_{\mathcal{A}_{\bf d}}^{\{r \}}$
for $r \leq s$ will all have the expected dimension by Theorem
\ref{partitionable}.   

This combinatorial technique for proving that secant varieties to certain
Segre varieties have the expected dimension was introduced by
Catalisano, Geramita and Gimigliano in \cite{Catalisano2002}.
As pointed out in \cite{Catalisano2002}, finding an $s$-element index set $I$ with pairwise Hamming distance greater than $2$ is equivalent to finding a placement of  $s$ rooks on a $(d_1 + 1) \times \cdots \times (d_n + 1)$ chessboard with the property that no two rooks attack each other or attack the same square on the board.  Our approach via triangulations  can be used to get information about
 further invariants (beyond dimension) of such secant varieties.
  \qed
\end{ex}

To conclude this section, we explore the existence of delightful triangulations for the class of rational normal scrolls.  While all the secant ideals in question are known to have nice determinantal presentations, not every scroll has a delightful term order.  This is somewhat surprising, considering our results on delightful term orders for minors and Pfaffians in Section \ref{sec:initial}.

Let  $\,{\bf \lambda} \,= \,(\lambda_1, \ldots, \lambda_n)\,$
be any vector of positive integers.  The rational normal scroll 
$S(\lambda)$ is the toric variety given by the parametrization:
$$x_{ij} = s^j t_i, \quad i = 1, \ldots, n, \quad j = 0, \ldots, \lambda_i.$$

The corresponding vector configuration equals
$$\mathcal{A}_\lambda = \{ j {\bf e}_0 \oplus {\bf e}_i \, \, | \, \, 1 \leq i \leq n
\,\,\, \mbox{and} \, \, 0 \leq j \leq \lambda_i \}
\quad \subset \quad \nn^{n+1}. $$
The toric ideal corresponding to this parametrization is denoted $I_\lambda$.

A determinantal presentation for the secant ideals $I_{\lambda}^{\{r\}}$ is well known.
Namely, for each
$i$ and $r$ with $\lambda_i \geq r$ let $M^{i,r}$ denote the $(r+1)
\times (\lambda_i - r + 1)$ Hankel matrix
$$M^{i,r} = \begin{pmatrix}
x_{i0} & x_{i1} & \ldots & x_{i,\lambda_i-r} \\
x_{i1} & x_{i2} & \ldots & x_{i, \lambda_i - r + 1} \\
\vdots & \vdots & \ddots & \vdots \\
x_{ir} & x_{ir+1} & \ldots & x_{i \lambda_i}
\end{pmatrix}.$$
If $\lambda_i < r$ then $M^{i,r}$ denotes the empty $(r+1) \times 0 $
matrix. The concatenation
$$ \,M^r \,\,\,= \,\,\,\left( M^{1,r} | M^{2,r} | \cdots | M^{n,r} \right) $$
is a matrix with $r+1$ rows and 
$\,\sum_{i=1}^r \lambda_i - n (r-1) \,$ columns.

\begin{thm}[\cite{Catalano1996, DePoi1996}] \label{prop:scrolldet}
The secant ideal $I^{\{r \}}_\lambda$ is generated by the  $(r + 1) \times (r
+ 1)$ minors of the matrix $M^r$.
\end{thm}

Our first result shows that delightful scrolls are rare.

\begin{prop}\label{thm:scroll}
If the ideal $I_{\bf \lambda}$ has a delightful term order then there exists an integer $m$ such that $\lambda_i \in \{m, m+1, m+2, m+3\}$ for all $i$.
\end{prop}

\begin{proof}
We first reduce to the two dimensional case.  This reduction is possible
because a delightful triangulation of a polytope is delightful for any face,
and the quadrangle
$\conv(\ca_{\lambda_i, \lambda_j})$ appears as a face of the polytope $\conv(\ca_\lambda)$.

To analyze the two dimensional case, we first must understand the full triangulations of the sets $\ca_{\lambda}$.  Each of these triangulations is lexicographic.  The full triangulations 
of $\ca_{\lambda_1,\lambda_2}$
 correspond to certain bipartite graphs.  Namely, aside from the edges like $(x_{1,i}, x_{1,i+1})$ and $(x_{2,i}, x_{2,i+1})$, the remaining edges form a bipartite planar spanning tree in the complete bipartite graph $K_{\lambda_1 + 1, \lambda_2 + 1}$.  Planar means that there is no pair of edges $(x_{1,i}, x_{2,j})$ $(x_{1,k}, x_{2,l})$ with $i < k$ and $j > l$.   An example of such a triangulation and the associated bipartite planar spanning tree appear in 
Figure 3.

\begin{figure}[h]
\begin{center}
\includegraphics{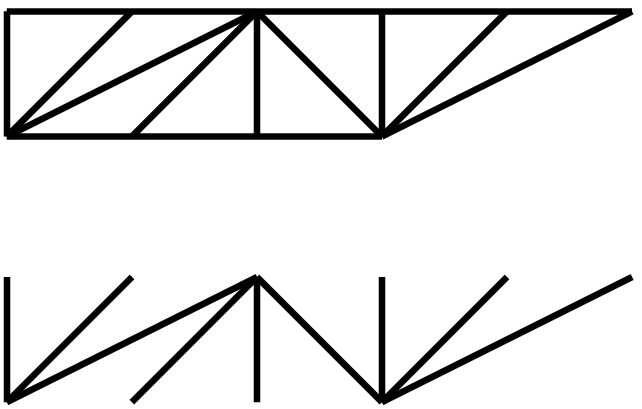}
\end{center}
{\bf Figure 3}. Triangulation for the scroll 
and the corresponding planar tree.
\end{figure}

If a triangulation of a scroll is delightful, then
the associated bipartite planar graph cannot possess certain induced subgraphs. 
We claim that, up to symmetry, these forbidden induced subgraphs are the following two:

\begin{enumerate}
\item  
The three edges  $(x_{1,i}, x_{2,j})$, $(x_{1,i}, x_{2,j+1})$, $(x_{1,i}, x_{2,j+2})$ 
where $1\leq i \leq \lambda_1 - 1$ 
and $ 0 \leq j \leq \lambda_2-2$.
\item The four edges
$(x_{1,0}, x_{2,0})$, $(x_{1,0}, x_{2,1})$, $(x_{1,0}, x_{2,2})$, $(x_{1,0}, x_{2,3})$. 
\end{enumerate}

Thus (1) is a claw $K_{1,3}$ which is not adjacent to a vertical boundary,
and (2) is a claw $K_{1,4}$ which is adjacent to a vertical boundary. 
 Note that the triangulation in
 Figure 3 contains both of these forbidden subgraphs.
If the graph of a triangulation contains the $K_{1,3}$ in case (1) then the ideal
$({\rm in }_\prec (I_\lambda))^{\{2\}}$ contains the monomial $x_{1,i-1} x_{1,i+1} x_{2,j+1}$.  However, by virtue of the fact that the full triangulation can be chosen to be lexicographic with $x_{1,i-1}$, $x_{1,i}$ and $x_{1,i+1}$ smaller than any other $x_{1,j}$, and appealing to Proposition \ref{prop:scrolldet}, we see that this monomial cannot be the leading monomial of any polynomial in $I_{\lambda}^{\{2\}}$.  A similar argument rules out the subgraph
$K_{1,4}$ in case (2).

To finish the proof, note that if $\lambda_1 < \lambda_2$ and $\lambda_1 +3 < \lambda_2$, then the induced graph of any full triangulation of $\ca_\lambda$ must contain one of the two forbidden subgraphs (or a subgraph symmetrically equivalent).
\end{proof}

We can, however, show the existence of a
delightful term order in the  special case
when all the $\lambda_i$ are equal.

\begin{thm}\label{lemma:scrollsame}
Suppose that $\lambda_1 = \lambda_2 = \cdots = \lambda_n$.
 Let $\prec$ be the lexicographic term order such that $x_{ij} \succ x_{kl}$ if $j < l$ or $j = l$ and $i < k$.  Then $\prec$ is delightful for $I_{\lambda}$.\end{thm}

\begin{proof}
The edge graph of every full triangulation of a configuration
$\mathcal{A}_\lambda$ is a chordal graph. This can be proved
by induction on $\sum_{i=1}^n \lambda_i$. Let $G_\lambda$ be
the complementary graph to that chordal graph. The initial
ideal ${\rm in}_\prec(I_\lambda) $ equals the edge ideal $I(G_\lambda)$.
Since chordal graphs are perfect, and the complements of perfect graphs
are perfect, it follows that $G_\lambda$ is a perfect graph.

For the particular lexicographic term order we have chosen,
the edges in the graph $G_\lambda$ are the pairs of the form $(x_{ij}, x_{kl})$ such that $j+1 < l$ or $j+1 = l$ and $i < k$.  It is the simplicity of the graph $G_\lambda$ which depends on $\lambda_i = \lambda_j$ for all $i$ and $j$.   To show that $\prec$ is delightful, we must show that for each clique of size $r$ in $G_\lambda$ there is a polynomial in $I_\lambda^{\{r-1\}}$ which has the clique as a leading term. Let $x_{i_1j_1} x_{i_2j_2} \cdots x_{i_r j_r}$ be such clique.  We may suppose that $j_1 < j_2 < \cdots < j_r$.  Consider the $r \times r $ matrix

$$M = 
\begin{pmatrix}
x_{i_1j_1}  & x_{i_2 j_2 -1} & x_{i_3 j_3 - 2} & \cdots & x_{i_r j_r - r + 1} \\
x_{i_1 j_1 + 1} & x_{i_2 j_2} & x_{i_3 j_3 - 1} & \cdots & x_{i_r j_r - r + 2} \\
x_{i_1 j_1 + 2} & x_{i_2 j_2 + 1} & x_{i_3 j_3} & \cdots & x_{i_r j_r - r + 3} \\
\vdots    &   \vdots   & \vdots  & \ddots  & \vdots \\
x_{i_1 j_1 + r - 1} & x_{i_2 j_2 + r - 2} & x_{i_3 j_3 + r -3} & \cdots &  x_{i_r j_r}
\end{pmatrix}.
$$

By construction, the polynomial $f = \det M$ belongs to $I_\lambda^{\{r-1\}}$ since it is one of the minors appearing in Proposition \ref{prop:scrolldet}.  Furthermore, $f$ is not identically zero.  The structure of $M$ implies that $f$ could be identically zero only if there were two identical columns in $M$.  This implies that there are indices $s$ and $t$ with $s < t$ such that $i_s = i_t$ and $j_s = j_t + t - s$.  But the conditions on the edges of $G_\lambda$ make it impossible for there to be a lexicographically ordered clique in $G_\lambda$ with these properties.  Furthermore, each indeterminate appearing in the matrix $M$ is a valid indeterminate in our polynomial ring.  This follows because all the second indices, $j_k \pm l$, lie between $j_1$ and $j_r$.  Finally, the term order $\prec$ selects the main diagonal as the leading term.  Hence $\prec$ is delightful for $I_\lambda$.
\end{proof}

\begin{cor}\label{cor:same1}
Suppose that there exists an integer $m$ such that $\lambda_i \in \{m,m+1\}$ for all $i$.  Then $I_\lambda$ has a delightful term order.
\end{cor}

\begin{proof}
The ideal $I_{\lambda}$ can be realized as the elimination ideal of $I_{\lambda'}$ where $\lambda'_i = m+1$ for all $i$.  The lexicographic ordering from Theorem \ref{lemma:scrollsame} realizes this elimination and thus the delightful property passes to this elimination ideal.
\end{proof}

In general, we do not know whether the converse to Proposition \ref{thm:scroll} is true.  However, we can show that it holds in the two dimensional case.

\begin{prop}
Suppose that $\lambda_2 \leq \lambda_1 \leq \lambda_2 +3$.  Then the ideal $I_{\lambda_1, \lambda_2}$ has a delightful term order.
\end{prop}

\begin{proof}
We will prove the case $\lambda_1 = \lambda_2 + 3$.  The case of $\lambda_1 = \lambda_2 + 2$ follows by the elimination argument used in the proof of Corollary \ref{cor:same1}, and the other two cases are proved in Theorem \ref{lemma:scrollsame} and Corollary \ref{cor:same1}.

Now introduce the lexicographic term order $\prec$, given by the rule
$$x_{10} \succ x_{11} \succ x_{20} \succ x_{12} \succ x_{21} \succ x_{13} \succ \cdots \succ x_{1j} \succ x_{2,j-1} \succ x_{1,j+1} \succ \cdots $$  $$\cdots \succ x_{1,\lambda+1} \succ x_{2,\lambda} \succ x_{1,\lambda+2} \succ x_{1,\lambda+3}.$$
We claim that this lexicographic term order is delightful for $I_{\lambda_1, \lambda_2}$.  To see this, let $x_{i_1j_1} , \ldots, x_{i_r j_r}$ be an independent set in the triangulation corresponding to this term order, arranged in decreasing lexicographic order.  Since every full triangulation of $\mathcal{A}_{\lambda_1,\lambda_2}$ is chordal and hence perfect, we must show that this independent set yields the initial term of a polynomial in $I^{\{r-1\}}_{\lambda_1, \lambda_2}$.  For a general pair of sequential elements in our independent set $x_{i_kj_k}, x_{i_{k+1}j_{k+1}}$, this happens if and only if either: $i_k = i_{k+1}$ and $j_k +1 < j_{k+1}$; or $i_k = 1, i_{k+1} = 2$, and $j_k \leq j_{k+1}$; or  $i_k = 2, i_{k+1} = 1$, and $j_k + 2 < j_{k+1}$.  The only exceptions to these rules come at the ends:  we cannot have the pairs $x_{10}, x_{20}$ or $x_{2\lambda_2}, x_{1\lambda_1}$ in an independent set. 

Now construct the matrix $M$ as in the proof of Theorem \ref{lemma:scrollsame}.  Our conditions on the sequence $x_{i_1j_1} , \ldots, x_{i_r j_r}$
 guarantee that all the entries in $M$ are valid indeterminates in our polynomial ring.  Furthermore, $f = \det M$ is not identically zero, and $f$ has leading term equal to $x_{i_1j_1}  \cdots x_{i_r j_r}$.  Thus, the lexicographic term order $\prec$ is delightful.
\end{proof}


\end{document}